\documentclass[a4paper, 12pt]{amsart}
\usepackage{amsmath, amssymb, eucal, amscd, amstext, enumerate}
\usepackage{amsfonts,amsthm}
\usepackage{mathrsfs}

\theoremstyle{plain}
\newtheorem{thm}{Theorem}[section]
\newtheorem{lem}[thm]{Lemma}
\newtheorem{eg}[thm]{Example}
\newtheorem{prop}[thm]{Proposition}
\newtheorem{cor}[thm]{Corollary}
\newtheorem{defn}[thm]{Definition}
\newtheorem{rem}[thm]{Remark}
\newtheorem{rem-ntn}[thm]{Remark and Notation}

\newenvironment{prf}{{\noindent \textbf{Proof:}\ }}{\hfill $\Box$\\ \smallskip}

\numberwithin{equation}{section}

\newcommand{\smnoind}{{\smallskip\noindent}}

\newcommand{\lsp}{{\rm span\ \!}}

\newcommand{\abs}[1]{\left\vert#1\right\vert}
\newcommand{\norm}[1]{\left\|#1\right\|}

\newcommand{\CP}{\mathcal{P}}

\newcommand{\CU}{\mathcal{U}}
\newcommand{\CI}{\mathcal{I}}

\newcommand{\CA}{\mathcal{A}}

\newcommand{\BR}{\mathbb{R}}

\newcommand{\BC}{\mathbb{C}}

\newcommand{\Spec}[1]{\mathrm{Sp}#1}

\begin{document}

\title[]{The number of Dirac-weighted eigenvalues of Sturm-Liouville equations with integrable potentials and an application to inverse problems}

\author{Xiao Chen \and Jiangang Qi$^\dag$}
\address{School of Mathematics and Statistics, Shandong University at Weihai \\Weihai 264209, P.R. China}
\email{chenxiao@sdu.edu.cn; qjgsduwh@126.com}\




\begin{abstract}

In this paper, we further Meirong Zhang, et al.'s work by computing the number of weighted eigenvalues for Sturm-Liouville equations, equipped with general integrable potentials and Dirac weights, under Dirichlet boundary condition. 
We show that, for a Sturm-Liouville equation with a general integrable potential, if its weight is a positive linear combination of $n$ Dirac Delta functions, then it has at most $n$ (may be less than $n$, or even be $0$) distinct real Dirichlet eigenvalues, or every complex number is a Dirichlet eigenvalue; in particular, under some sharp condition, the number of Dirichlet eigenvalues is exactly $n$. 
Our main method is to introduce the concepts of characteristics matrix and characteristics polynomial for Sturm-Liouville problem with Dirac weights, and put forward a general and direct algorithm used for computing eigenvalues.
As an application, a class of inverse Dirichelt problems for Sturm-Liouville equations involving single Dirac distribution weights is studied.

\vspace{04pt}

\vspace{04pt}

\end{abstract}

\begin{thanks} {$\dag$ the corresponding author, qjgsduwh@126.com}
\end{thanks}

\maketitle

\section{Introduction}
\medskip

The inverse spectrum problems of Sturm-Liouville equations have always been the subjects of intense scholarly research.
The classical inverse problems focus on the studies of how to determine the potential uniquely using suitable spectral data.
In 1946, G.\ Borg wrote the first important paper \cite{Bo46} on this subject.
From then on, a number of literatures have been devoted to this direction (cf.  \cite{CCPR87, FWW16, Ha76, Ho74, Hor05, KK97, PT87, QC16, Yu00} etc).

Many eigenvalue problems often involve Dirac delta functions, although the coefficients of the classical Sturm-Liouville equations are considered in the space $L^1[a,b]$. 
Even for regular Sturm-Liouville problems,  the corresponding optimization problem of weights or potentials were also frequently related to Dirac Delta functions (cf. \cite{GQ17, GQ20, Gu19, Hin13, MZ13, WZZ20, ZWMQX18} etc).
More to the point, in physics, it is very important and common when the coefficients of second-order ordinary differential equations contain Dirac distributions. In \cite{KP31}, physicists Kronig and Penney introduced an equation with Dirac delta potential for describing quantum mechanics in crystal lattices; in \cite{AN07}, Dirac Delta functions appeared in a continuum physical model of nonlocal point interactions. 
Moreover, since every integrable function can be weakly approximated by distributions, and the eigenvalues is continuous with respect to distributions or measures under $w^*-$topology (see \cite[Lemma 4.1]{GQ20} and \cite[Theorem 1.3]{MZ13}), the eigenvalue problems of differential equations can be transformed into the problems of finding the roots of algebraic equations. 
This provides an effective approach to the approximate evaluation of eigenvalues.

In the present paper, we study a class of general \emph{Dirac-weighted Sturm-Liouville eigenvalue problems}, in which the potentials are general integrable functions and the weights are linear combinations of finitely many Dirac Delta functions. 
This kind of equations has explicit and significant physical backgrounds, for instance,  a forced vibration system of a string attached with finitely many particles (see Section~\ref{sec:pre}).
Moreover, in \cite{GQ20}, the authors also described another similar model, i.e., a force-free  vibration system which made up of a string fixed at one end and finitely many masses connected to a spring at the other end. 
And then, a series of new inverse problems and new optimization problems related to our topic accordingly arise.  

In the remainder of this section, we give a brief introduction to the main content of our paper, which is contained in Section~ \ref{number-eigenv}, Section~\ref{sec:inverse-prob} and Section~\ref{sec:question}. Besides, Section~\ref{sec:pre} describes our problems in detail and defines their solutions.

\noindent
{\bf Part I. The number of Dirac-weighted eigenvalues of Dirichlet problems with general integrable potentials.}
In \cite{ZWMQX18}, Zhang and his collaborators, using the difference method, proved that, a Sturm-Liouville equation with zero potential must have exactly $n$ Dirichlet eigenvalues, if its weight is a linear combination of $n$ Dirac Delta functions. Unfortunately, the difference method used in \cite{ZWMQX18} is only applicable for the zero potential case, and the corresponding results about the number of Dirichlet eigenvalues does not always holds in general.

In Section~\ref{number-eigenv}, we will in a new way generalize Zhang's results to the general case when the potential is a general integrable function (see Theorem~\ref{thm:num-eigen}). This generalization is highly non-trivial. 
Our results reveal more essences of this eigenvalue problem.
More precisely, we prove that, for any Sturm-Liouville equation with a non-zero integrable potential,  if the weight is a positive linear combination of $n$ Dirac Delta functions, then the set of its Dirichlet eigenvalues either consists of at most $n$  (may be less than $n$, or even be $0$) different real numbers, or is the whole complex plane (see Theorem~\ref{thm:num-eigen}(a)(e)). 
And the number of its Dirichlet eigenvalues strictly depends on some key hypotheses. 
In particular, when one of these hypotheses holds, the number of Dirichlet eigenvalues is exactly $n$ (see Theorem~\ref{thm:num-eigen}(c)(d)).
It is worth to notice that, in \cite[Lemma 3.3]{ZWMQX18}, the absence of the potential happens to guarantee that the key hypotheses automatically hold, and makes the direct calculation on eigenvalues easy (see Section~\ref{subsubsec:zero-q}).

In addition, we establish a systemical method for dealing with Sturm-Liouville eigenvalue problems involving Dirac weights, which most probably produces a new numerical method of computing eigenvalues (see Remark~\ref{rem:alg}).
Simply put, by tracking solutions of our equations for given initial values, we construct characteristics  symmetric matrices and characteristics polynomials with respect to such eigenvalue problems above, and thereby derive that these Sturm-Liouville Dirac-weighted eigenvalue problems are essentially equivalent to the problems of computing roots of their characteristics polynomials (equivalently, eigenvalue problems of their characteristics matrices), see Theorem~\ref{thm:eigen-E-mat-poly}. 
In the special case when the potential is zero, our approach just coincides with the difference method employed in \cite{ZWMQX18} (see Section~\ref{subsubsec:zero-q}). 

\noindent
{\bf Part II. Inverse Dirichelt problems for SLPs involving single Dirac distribution weights.} 
As a direct consequence of Section~\ref{number-eigenv}, we can see that, for a Sturm-Liouville equation on $[0,1]$, under some suitable conditions, if its weight is a single Dirac Delta function at some point $t\in(0,1)$,  then there exists a unique strictly positive Dirichlet eigenvalue, say $\lambda(t)$. Here a set $\{\lambda(t)\}_{t\in(0,1)}$ are called \emph{complete spectrum data}.
See Section~\ref{subsubsec:singelDirac}, Lemma~\ref{lem:property-eigen} and Corollary~\ref{cor:eigen-lim}.

In Section~\ref{sec:inverse-prob}, as an application, we solves a class of inverse problems: the potential can be uniquely recovered by the complete spectrum data set $\{\lambda(t)\}_{t\in(0,1)}$ mentioned above. 
In fact, this result can be further strengthened. 
For any suitable function $f(t)$ on $(0,1)$, called a \emph{spectrum-like function} (see Definition~\ref{defn:spec-like}), we can uniquely construct a Sturm-Liouville equation involving a single Dirac distribution weight such that all its Dirichlet eigenvalues over $(0,1)$ are given by $\{f(t)\}_{t\in(0,1)}$. 
In fact, we have owned an equivalence between spectrum-like functions and complete spectrum data (Theorem~\ref{thm:q-spec-like}).  
Also, we give a simple example for describing the physical meaning of this equivalence above (see Example~\ref{eg:vibration}).
Finally, it is worth mentioning that our work will certainly lead to a batch of new inverse problems and new optimization problems, some of which are listed in Section~\ref{sec:question}.

\bigskip
\section{Problem statement and preliminaries}\label{sec:pre}
\medskip

Let us begin with a forced vibration model, which is made up of a string and a particle $P$ attached on this string, and executes small transverse vibrations about the position of stable equilibrium, the interval $0\leqslant x\leqslant l$ of the $x$-axis.  
The motion is said to be ``small" in the sense that higher powers of the function $u(x,t)$ below and of its derivatives may be neglected compared with lower powers (cf. \cite[Section IV.10.2]{CH89}). 
Besides, we also assume that the mass of the string, compared with the mass $m$ of the particle, is negligible and set as $0$. And let $x=\xi$ be the coordinate of the particle $P$.

Denote by $p(x)$ the modulus of elasticity of the string multiplied by the cross-sectional areas, by $u(x,t)$ the perpendicular deflection of a point on the string from equilibrium position at time $t$, and by $F(x,t)$ the external force. 

First, let us consider the case when the string is fixed at its end points, i.e., $u(0,t)=u(l,t)=0$, which is called the boundary condition of the system here.  
At any time $t$, the kinetic energy of the string is given by the integral 
$$T(t)=\frac{1}{2}mu^2_t(\xi,t)=\frac{1}{2}\int_0^lm\delta(x-\xi)u^2_t(x,t)\, dx,$$ 
where $\delta(x-\xi)$ is the Dirac delta function at $\xi$ defined as in \eqref{defn:delta} below, and the potential energy takes the form 
$$U(t)=\frac{1}{2}\int_0^l \{p(x)u^2_x(x,t)+F(x,t)u(x,t)\}\, dx.$$

Next, we consider a special external force determined by a non-homogeneous intensity $q(x)$ along $x$-axis and the perpendicular deflection $u(x,t)$, i.e., $F(x,t)=q(x)u(x,t)$,  where $q(x)$, here called the  stationary  potential of the non-uniform force field $F$, is an integrable function on $[0,l]$. So, the total energy is 
$$E(t)=T(t)-U(t)=\frac{1}{2}\int_0^l \{m\delta(x-\xi)u^2_t(x,t)-p(x)u^2_x(x,t)-F(x,t)u(x,t)\}\, dx.$$

Because of the law of conservation of energy, we know $E'(t)=0$. Then, from the standard calculus, we obtain 
$$0=E'(t)=\int_0^l u_t(x,t)\{m\delta(x-\xi)u_{tt}(x,t)+(p(x)u_x(x,t))_x-q(x)u(x,t)\}\, dx,$$
which gives the partial differential equation of the vibrating string
\begin{equation}\label{eqn:vib-mass}
-(pu_x)_x(x,t)+q(x)u(x,t)=m\delta(x-\xi)u_{tt}(x,t).
\end{equation}
Separating variables in \eqref{eqn:vib-mass}, we write the solution $u$ in the form $v(x)e^{i\omega t}$ (cf. \cite[Section V.3]{CH89}), where $\omega$ is the vibrating frequency. For the amplitude $v(x)$, setting $\lambda=\omega^2$, from \eqref{eqn:vib-mass}, we then have the ordinary differential equation with Dirichlet boundary condition
\begin{equation}\label{eqn:vib-mass-v}
-(pv')'(x)+q(x)v(x)=\lambda m\delta(x-\xi)v(x)\text{ on }[0,1],\quad v(0)=v(l)=0.
\end{equation}
It is a Sturm-Liouville equation with a single Dirac weight and a integral potential, under Dirichlet boundary condition. 

More generally, a forced vibration system of a string attached with finitely many particles and two fixed ends can be characterized as the problem {\bf(E)} below.  And, due to Liouville transformation (see e.g.\cite[Section V.3]{CH89}), we may as well assume that $p(x)\equiv1$. 

In the present paper, we mainly consider a class of  Sturm-Liouville problems (SLPs) involving Dirac weights and integrable potentials as follows:
\begin{equation}\label{defn:E}
   \quad\quad -y'' + qy = \lambda \delta_{n,\vec{m},\vec{t}}y,\ y=y(x),\ x\in[0,1]
\end{equation}
associated to the Dirichlet boundary condition
\begin{equation}\label{defn:Eb}
y(0)=0=y(1)
\end{equation}
by the definitions of solutions of SLPs.
Here $q\in L^1([0,1],\BR)$, and set 
\begin{equation}\label{defn:dirac-weight}
\delta_{n, \vec{m},\vec{t}}(x) := \sum_{i=1}^{n} m_i \delta(x-t_i),
\end{equation}
where 
\begin{equation}\label{defn:dirac-contants}
\begin{array}{l}
   \vec{m} := (m_1, \cdots, m_n), \ m_i> 0\ (i=1,2,\cdots, n), \\
  \vec{t} := (t_0, t_1, \cdots t_n, t_{n+1}), \ 0 = t_0 < t_1 < \cdots < t_n < t_{n+1} = 1,
\end{array}
\end{equation}
and $\delta(x-t_i)$ is the canonical Dirac Delta function at the point $t_i$. 

\noindent
{\bf Notation:} denote by {(\bf E)} the SLP \eqref{defn:E}-\eqref{defn:Eb}.

First of all, similar to \cite[Section~2]{GQ20}, we give the definition of solutions of the SLPs with distributions.
Set
\begin{equation}\label{defn:D}
D[0,1]:=\lsp_{\BR}\Big\{L^1([0,1],\BR)\cup\big\{\delta(x-t_j):\
t_j\in[0,1],\ \mathbb N\ni j<\infty\big\}\Big\}, 
\end{equation}
where $\BR$ is the field of real numbers, $L^1([0,1],\BR)$ is the space of Lebesgue integrable, real valued functions on $[0,1]$,
and $\delta(x-t)$ is the Dirac delta function at point $t\in[0,1]$
defined by
\begin{equation}\label{defn:delta}
\delta(x-t)=
\begin{cases}
\infty, &x=t,\\
          0, &x\neq t
\end{cases}
\ \text{and}\ \int_I\delta(x-t)\, dx=1,\ \ \forall I\subset[0,1]
\ \text{and}\ t\in I,
\end{equation}
which can also be considered as the Radon-Nikodym derivative
of Dirac measure
\begin{equation*}
\delta_{t}(x)=
\begin{cases}
0, &x\in[0,t),\\
          1, &x\in[t,1].
\end{cases}
\end{equation*}
Then for $w\in D[0,1]$, which is called a {\bf\emph{distribution (or
distribution density)}} on $[0,1]$,
\begin{equation}\label{defn:D-weight}
w(x)=r\rho(x)+\sum_{j=1}^{n}r_j\delta(x-t_j),
\end{equation}
where $r, r_j$ are constants,
$\rho\in L^1[0,1]$ and $0\leqslant t_1<t_2<\cdots<t_n \leqslant1$.

Define the norm of $w$ by
\begin{equation}\label{defn:D-norm}
\|w\|_D=|r|\int_0^1|\rho|\mathrm dx+\sum_{j=1}^n|r_j|.
\end{equation}

Let $D_0[0,1]$ be a subset of $D[0,1]$ consisting of the distributions, called {\bf\emph{$D_0$-distributions}}, like
\begin{equation}\label{eq:D0-weight}
r\rho(x)+\sum_{j=1}^{n}r_j\delta(x-t_j), \text{ i.e., }r\rho(x)+\delta_{n,\vec{r},\vec{t}}
\end{equation}
where $r, r_j>0,$
$\rho\in L^1[0,1]$ and $0<t_1<t_2<\cdots<t_n<1.$
It is clear that $D_0[0,1]$ excludes any linear combination of $\delta(x-0)$ and $\delta(x-1)$. 
The physical reason for considering $D_0[0,1]$ instead of  $D[0,1]$ is that it makes little sense to put any mass at two fixed ends. And from the point view of mathematics,  Example 3.1 in \cite{ZWMQX18} told us that both of $\delta(x-0)$ and $\delta(x-1)$ are of no account to any result we want.

Let $w\in D_0[0,1]$ with the form of \eqref{eq:D0-weight}, $q\in L^1[0,1]$ and $\ I_i:=(t_i,t_{i+1})\ (i=0,1,\cdots,n).$

{\bf \emph{A solution of the SLP with a $D_0$-distribution $w$}}
\begin{equation}\label{defn:E-D}
\quad\quad -y'' + qy = \lambda wy,\ y=y(x),\ x\in[0,1],
\end{equation}
associated to the boundary condition \eqref{defn:Eb}
is the function $y$ satisfying
$$
y\in\left\{z: z\in AC[0,1], \ z'\in \bigcup^{n}_{i=0} AC(I_i),
\ \exists z'(t_j^\pm),\ j=1,\cdots,n\right\}
$$
and
\begin{equation}\label{defn:E-D-bdd}
\left\{\aligned
&-y''(x)+q(x)y(x)=\lambda r\rho(x)y(x),\ x\not=t_j, x\in[0,1],\\
&y'(t_j^-)-y'(t_j^+)=\lambda r_j y(t_j), \ 1\leqslant j\leqslant n,\\
&y(0)=0=y(1),\\
\endaligned\right.
\end{equation}
where all $y'(t_j^+)$ and $z'(t_j^+)$ (resp., $y'(t_j^-)$ and $z'(t_j^-)$) denote the classical right-derivative (resp., left-derivative), 
$AC(I)$ denotes the space of absolutely continuous real-valued functions
on an subinterval $I\subset[0,1]$, and $y'(0^-)$ (resp., $y'(1^+)$) is written as $y'(0)$ (resp., $y'(1)$).

\noindent
{\bf Notation:} denote by ${(\bf E_{D})}$ the SLP \eqref{defn:E-D}-\eqref{defn:E-D-bdd}.

Obviously, the Dirichlet problem {\bf (E)} is actually a special case of the above SLP ${(\bf E_{D})}$ by setting 
\begin{equation}\label{defn:dirac-setting}
\rho=0,\ r_j=m_j\ (j=1,2,\cdots,n),
\end{equation}
where all $t_i\ (i=0,1,\cdots,n+1)$ and $m_j\ (j=1,2,\cdots,n)$ are defined as \eqref{defn:dirac-contants}. The number $\lambda\in\BC$ is an eigenvalue of {\bf (E)} 
if and only if the problem has a nontrivial
solution, $E(x,\lambda)$, which is called an eigenfunction corresponding to $\lambda$. 
Consequently, we have grounds to give the following definition of a solution of the SLP {\bf (E)}.

\begin{defn}\label{defn:solution}
{\bf \emph{A solution of the SLP {\bf (E)}}} is defined as a solution of ${(\bf E_{D})}$ with $w$ being $\delta_{n,\vec{m},\vec{t}}$ given in \eqref{defn:dirac-setting}.
\end{defn}

\begin{rem}\label{rem:mes-form}
For any $q\in L^1[0,1]$, its distribution, namely,
$\mu_q(x):=\int_0^x q(t)\, dt,$ for any $x\in[0,1],$
defines an absolutely continuous measure with respect to lesbegue measure on $[0,1]$.
And denote by $\Delta_{n,\vec{m},\vec{t}}$ the purely discontinuous measure (also known as completely singular measure) $\sum_{j=1}^{n} r_j \delta_{t_j}$ induced by $\delta_{n,\vec{m},\vec{t}}$.
Then, we can see that, the equation \eqref{defn:E-D} is in fact the measure differential equation (MDE) as follows:
\begin{equation}\label{eqn:mes-form}
-dy^{\bullet}(x)+y(x)\, d\mu_q(x)=\lambda y(x)\, d(r\mu_\rho+\Delta_{n,\vec{m},\vec{t}})(x),
\end{equation}
where $y$ is defined as in \cite[(1.3)-(1.4)]{ZWMQX18}, and $y^{\bullet}$ stands for the generalized right derivative of $y$. For more information on measures and eigenvalue theory of second-order MDEs, see \cite{Fol99}, \cite{HR94}, \cite{MZ13} and \cite{ZWMQX18}. But, for simplicity, we still continue to use the previous notations and terminology.
\end{rem}

\bigskip
\section{The number of Dirac-weighted eigenvalues of Dirichlet problems with integrable potentials}\label{number-eigenv}

In this section, we will study the structure of Dirichlet Dirac-weighted eigenvalues of the problem {\bf (E)}.

\bigskip
\subsection{Characteristics polynomials and characteristics matrices}

Set two problems as follows:
\begin{equation*}
  \begin{array}{l}
  {\bf (E_{\left[\xi,\eta \right]})}:\left\{ \begin{array}{ll}
    -y'' + qy = 0 & \text{ defined on } \left[\xi, \eta \right] \subsetneq [0, 1]\\
    y(\xi) = y(\eta) = 0
  \end{array}  \right.; \\
  {\bf (E_0)}:\left\{ \begin{array}{ll}
    -y'' + qy = 0 & \text{ defined on } [0,1]\\
    y(0)=y(1) = 0
  \end{array}  \right..
  \end{array}
  \end{equation*}
Denote by $\Spec (E_*)$ the set of eigenvalues of certain eigenvalue problem $E_*$, and by $\#\Spec (E_*)$ the number of the eigenvalues.

Let $\phi,\psi$ be the two linearly independent solutions of $-y'' + qy = 0$ on $\left[0,1\right]$ satisfying $\phi(0)=0$, and $W[\phi,\psi]$ be their Wronskian determinant, which must be non-zero, say $\omega \neq 0$.
Let $E(x,\lambda)$ be the eigenfunction of {\bf (E)} corresponding to  the eigenvalue $\lambda$. 

For the problem {\bf (E)}, denote $E_i(x,\lambda):= E(x,\lambda)|_{\left[t_i,t_{i+1}\right]}(i = 0,1,2,\cdots, n)$, which is actually a solution of $-y'' + qy = 0$ on $[t_i,t_{i+1}]$, can be presented as is a linear combination of $\phi$ and $\psi$, i.e. $E_i(x,\lambda) = \alpha_i \phi(x) + \beta_i \psi(x)$ for some $\alpha_i, \beta_i\in\BR$. From the uniqueness of solutions for initival value problem, we easily have the following result.

\begin{lem}\label{lem:dim-eigenspc}
$\dim V_{\lambda} = 1$ where $V_{\lambda}$ be the eigenspace of {\bf (E)} corresponding to $\lambda$. 
\end{lem}

At every nodal point $t_i$ $(i=1,\cdots,n)$, by the continuity and differentiability of eigenfunctions, we have the following basic relation:
\begin{equation*}
  \text(\bf R):
  \left\{ \begin{array}{ll}
    \alpha_i \phi(t_i) + \beta_i \psi(t_i) = \alpha_{i-1} \phi(t_i) + \beta_{i-1} \psi(t_i) \\
    \alpha_i \phi'(t_i) + \beta_i \psi'(t_i) = \alpha_{i-1} \phi'(t_i) + \beta_{i-1} \psi'(t_i)
    - \lambda m_i \left[\alpha_{i-1} \phi(t_i) + \beta_{i-1} \psi(t_i) \right] 
  \end{array}  \right.,
  \end{equation*}
   i.e., 
   $$\begin{pmatrix}
  \alpha_i \\
  \beta_i
\end{pmatrix} \overset{\bf (R_1)}{=} M_i(\lambda,\phi,\psi)\begin{pmatrix}
  \alpha_{i-1} \\
  \beta_{i-1}
\end{pmatrix}\overset{\bf (R_2)}{=} \begin{pmatrix}
  \alpha_{i-1} \\
  \beta_{i-1}
\end{pmatrix} + \frac{\lambda}{\omega} m_i N_i(\lambda,\phi,\psi) \begin{pmatrix}
  \alpha_{i-1} \\
  \beta_{i-1}
\end{pmatrix} ,\quad i = 1, \cdots , n+1, $$
where 
$$M_i(\lambda,\phi,\psi) = \left(\begin{matrix}
  1 + \frac{\lambda}{\omega} m_i (\phi \psi)(t_i) & \frac{\lambda}{\omega} m_i \psi^2(t_i) \\
  -\frac{\lambda}{\omega} m_i \phi^2(t_i) & 1 - \frac{\lambda}{\omega} m_i (\phi \psi)(t_i)
\end{matrix}\right),$$ 
and
$$N_i(\lambda,\phi,\psi) = \left(\begin{matrix}
  (\phi \psi)(t_i) & \psi^2 (t_i) \\
  -\phi^2(t_i) & -(\phi \psi)(t_i) \\
\end{matrix}\right) = \left(\begin{matrix}
  \psi(t_i) \\
  -\phi(t_i)
\end{matrix}\right) \left(\begin{matrix}\phi(t_i) & \psi(t_i)\end{matrix}\right).$$

\noindent
{\bf Fact:} \emph{$M_i \in SL_2(\BR)$, which implies that every $(\alpha_i, \beta_i)$$(i=1,2,\cdots, n)$ is a non-zero vector unless $(\alpha_0, \beta_0)=(0,0)$.}

Here $\alpha_i$ (resp. $\beta_i$) should be precisely written as $\alpha_i(\phi,\psi,n,\vec{m},\vec{t})$ (resp. $\beta_i(\phi,\psi,n,\vec{m},\vec{t})$), because these two coefficients depend on $\phi$, $\psi,n$, $\vec{m}$ and $\vec{t}$ set above. For convenience, we still simply write these two coefficients as $\alpha_i$ and $\beta_i$, if there is no ambiguity.

Define 
\begin{equation}\label{eq:Discrim}
D_{\eta, \xi}(\phi,\psi) := \det\left(\begin{matrix}
  \phi(\eta) & \phi(\xi) \\
  \psi(\eta) & \psi(\xi) \\
\end{matrix} \right)=\left(\begin{matrix}\phi(\eta) & \psi(\eta)\end{matrix}\right)
\left(\begin{matrix}
  \psi(\xi) \\
  -\phi(\xi)
\end{matrix}\right)
\end{equation}
for any $0 \le \xi < \eta \le 1$.  Here  we call this determinant {\bf\emph{discriminant}} for ${\bf (E_{\left[\xi, \eta\right]})}$ with respect to $\phi$ and $\psi$, since it is well known from the basic knowledge of ordinary differential equations that $D_{\eta, \xi} \ne 0$ if and only if ${\bf (E_{\left[\xi, \eta\right]})}$ only has the zero solution.  

For convenience, let $M_i:=M_i(\lambda,\phi,\psi)$, $N_i:=N_i(\lambda,\phi,\psi)$ and $D_{\eta, \xi}:=D_{\eta, \xi}(\phi,\psi)$, if there is no ambiguity. 

By the direct computation, we get the following useful lemma.
\begin{lem}\label{lem:Discrim-relation}
For any $0\leqslant \eta<a<b<\xi\leqslant 1$, one has 
$D_{\xi,a}\cdot D_{b,\eta}-D_{\xi,b}\cdot D_{a,\eta}=D_{\xi,\eta}\cdot D_{b,a}.$
\end{lem}

For any $1\leqslant K\leqslant n$, after $K$th iteration of  the relation ${\bf (R_2)}$, applying \eqref{eq:Discrim}, we have, by induction, that

\begin{equation}\label{eq:alpha-beta-K}
\begin{aligned}
  \left(\begin{matrix}
    \alpha_K \\
    \beta_K
  \end{matrix}\right) =& \left[I + \frac{\lambda}{\omega} \sum_{i = 1}^{K}m_i N_i + \sum_{l = 2}^{K-1} \left(\frac{\lambda}{\omega}\right)^l \sum_{1 \le i_1 < \cdots < i_l \le K} \prod_{k=1}^l m_{i_k} \prod_{k=l}^{1}N_{i_k} + \left(\frac{\lambda}{\omega}\right)^K \prod_{k=1}^{K}m_k \prod_{k=K}^{1}N_k\right] \left(\begin{matrix}
    \alpha_0 \\
    \beta_0
  \end{matrix}\right) \\
  =& \left[I + \frac{\lambda}{\omega} \sum_{i=1}^{K}m_i N_i + \sum_{l = 2}^{K-1} \left(\frac{\lambda}{\omega}\right)^l \sum_{1 \le i_1 < \cdots < i_l \le K} \prod_{k=1}^l m_{i_k} \prod_{k=2}^{l}D_{t_{i_k}, t_{i_{k-1}}} \cdot \left(\begin{matrix}
    \psi(t_{i_l}) \\
    -\phi(t_{i_l})
  \end{matrix}\right) \left(\begin{matrix}\phi(t_{i_1}) & \psi(t_{i_1})\end{matrix}\right) \right.\\
  &+ \left. \left(\frac{\lambda}{\omega}\right)^K \prod_{i=1}^{K}m_k \prod_{i=2}^{K}D_{t_i, t_{i-1}} \left(\begin{matrix}
    \psi(t_K) \\
    -\phi(t_K)
  \end{matrix}\right) \left(\begin{matrix}\phi(t_1) & \psi(t_1)\end{matrix}\right) \right] \left(\begin{matrix}
    \alpha_0 \\
    \beta_0
  \end{matrix}\right) ,
\end{aligned}
\end{equation}
where $I = \left(\begin{matrix}
  1 & 0 \\
  0 & 1 \\
\end{matrix}\right)$.

If choosing $\begin{pmatrix}\alpha_0 \\ \beta_0 \end{pmatrix} = \begin{pmatrix} 1 \\ 0 \end{pmatrix}$, by \eqref{eq:alpha-beta-K}, we have, for any $1\leqslant K\leqslant n$,
\begin{equation} \label{eq:alpha-K}
  \begin{aligned}
    \alpha_K =& \,\,
    1 + \frac{\lambda}{\omega} \sum \limits_{k=1}^{K} m_k (\phi \psi) (t_k) + \sum \limits_{l=2}^{K-1}\left(\frac{\lambda}{\omega}\right)^l \sum \limits_{1 \le i_1 < \cdots < i_l \le K} \prod \limits_{k=1}^{l} m_{i_k} \prod \limits_{k=2}^{l} D_{t_{i_k}, t_{i_{k-1}}} \psi(t_{i_l})\phi(t_{i_1}) \\
    &+ \left(\frac{\lambda}{\omega}\right)^K \prod \limits_{k=1}^K m_k \prod \limits_{k=2}^{K} D_{t_k, t_{k-1}} \psi(t_K) \phi(t_1), \\
  \end{aligned}
\end{equation}     
    
\begin{equation} \label{eq:beta-K}
  \begin{aligned}    
    \beta_K =&
    -\frac{\lambda}{\omega} \sum_{k=1}^{K} m_k \phi^2(t_k) - \sum_{l=2}^{K-1} \left(\frac{\lambda}{\omega}\right)^l \sum_{1 \le i_1 < \cdots < i_l \le K} \prod_{k=1}^{l} m_{i_k} \prod_{k=2}^{l} D_{t_{i_k}, t_{i_{k-1}}} \phi(t_{i_l})\phi(t_{i_1}) \\
    &- \left(\frac{\lambda}{\omega}\right)^K \prod_{k=1}^{K} m_k \prod_{k=2}^{K} D_{t_k, t_{k-1}} \phi(t_K) \phi(t_1),
  \end{aligned}
\end{equation} 
and
\begin{equation}\label{eqn:eigenfunction}
E(x, \lambda) = \left\{\begin{array}{ll}
  \phi(x), & x \in [0,t_1]\\
  \alpha_i \phi(x) + \beta_i \psi(x), & x \in [t_i, t_{i+1}],\ i = 1,2, \cdots, n
\end{array} \right..
\end{equation}

Define
\begin{equation}\label{defn:charpoly}
p(\lambda,\phi,\psi) :=\alpha_n \phi(1) + \beta_n \psi(1) = E(1,\lambda),
\end{equation}
which is a polynomial of degree at most $n$ with respect to the argument $\lambda$.

\begin{defn}\label{defn:charmat-charpoly}
The polynomial $p(\lambda,\phi,\psi)$ in \eqref{defn:charpoly} is called a {\bf\emph{characteristics polynomial}} of {\bf (E)}. If there exists a symmetric matrix $X(\phi,\psi)$ such that $p(\lambda,\phi,\psi)=c\det(X(\phi,\psi)-\lambda\cdot I)$ where $I$ is a unit matrix and $c$ is a non-zero constant, then $X(\phi,\psi)$ is called a {\bf\emph{characteristics matrix}} of {\bf (E)}.  For simplicity, we usually denote $p(\lambda,\phi,\psi)$ (resp., $X(\phi,\psi)$) by $p(\lambda)$ (resp., $X$) if both of $\phi$ and $\psi$ are fixed.
\end{defn}
Note that  the characteristics polynomial $p(\lambda)$ of {\bf (E)} is unique up to multiplication by scalars; and the characteristics matrix $X$ of  {\bf (E)} is unique up to similarity equivalence, if it exists. Clearly, due to Dirichlet boundary and \eqref{defn:charpoly}, we can see that the number $\lambda$ is a Dirichlet eigenvalue of {\bf (E)} if and only if $E(1,\lambda)=0$ for $E(x,\lambda)$ given in \eqref{eqn:eigenfunction}. So the following theorem is obtained.

\begin{thm}\label{thm:eigen-E-mat-poly}
The number $\lambda$ is a Dirichlet eigenvalue of {\bf (E)}, if and only if $\lambda$ is a root of the polynomial equation $p(\lambda)=0$, if and only if $\lambda$ is an eigenvalue of the symmetric matrix $X$ provided that its characteristics matrix $X$ exists.
\end{thm}

\begin{rem}\label{rem:alg}{\bf (Algorithm of calculating eigenvalues)} From the above arguments, for the problem {\bf(E)}, we can compute its Dirichlet eigenvalues by following algorithm:
\begin{itemize}
  \item[Step\,1:]  find and fix two linearly independent solutions $\phi$ and $\psi$ of $-y''+qy=0$ on $[0,1]$ satisfying $\phi(0)=0$.
  \item[Step\,2:]  for $\vec{t}$ in {\bf(E)}, work out  by \eqref{eq:Discrim} $\frac{n(n+1)}{2}$ discriminants $D_{t_j, t_i}$  $(1\leqslant i<j\leqslant n)$, and $\phi(t_k)$, $\psi(t_k)$ $(1\leqslant k\leqslant n)$.
  \item[Step\,3:]  obtain $p(\lambda)$ defined as \eqref{defn:charpoly} by calcultating both of $\alpha_n$ and $\beta_n$ in \eqref{eq:alpha-K} and \eqref{eq:beta-K}.
  \item[Step\,4:]  solve out the roots of $p(\lambda)=0$. The set of the eigenvalues of {\bf(E)} exactly consists of these roots.
\end{itemize}
\end{rem}

In the next section, we will use this algorithm above to compute the number of Dirichlet eigenvalues.

\bigskip
\subsection{The number of Dirichlet Dirac-weighted eigenvalues}\label{subsec:num-eigen}
At the beginning of this part, for any two fixed linearly independent solutions $\phi,\psi$ of $-y'' + qy = 0$ on $\left[0,1\right]$ and the partition $0 = t_0 < t_1 < \cdots < t_n < t_{n+1} = 1$ in {\bf (E)}, we propose three hypotheses:
\begin{itemize}
  \item[] ${\bf (H_0)}:$ the problem ${\bf (E_0)}$ only has the zero solution, i.e. $D_{1,0} \ne 0$.
  \item[] ${\bf (H)}:$ for all $i \in \{0, 1, \cdots n\}$, we have that all ${\bf (E_{[t_i, t_{i+1}]})}$ only have the zero solution, i.e. $D_{t_{i+1}, t_i} \ne 0$.
   \item[] ${\bf (H_1)}:$ for any $i \in \{1, \cdots n\}$, we have that both of ${\bf (E_{[t_0, t_i]})}$ and ${\bf (E_{[t_i, t_{n+1}]})}$ have the non-zero solution, i.e. $D_{1, t_i}=D_{t_i, 0}=0$.
\end{itemize}

\begin{rem}\label{rem:relation-hypo}
The hypothesis ${\bf (H_1)}$ holds, if and only if, for all $0\leqslant i<j \leqslant n+1$, all ${\bf (E_{[t_i, t_j]})}$ have the non-zero solution, i.e. $D_{t_j, t_i}=0$. Indeed, the hypothesis ${\bf (H_1)}$ is equivalent to the fact that $\frac{\phi(0)}{\psi(0)}=\frac{\phi(t_1)}{\psi(t_1)}=\frac{\phi(t_2)}{\psi(t_2)}=\cdots=\frac{\phi(t_n)}{\psi(t_n)}=\frac{\phi(1)}{\psi(1)}$.  Clearly, in this case, neither {\bf (H)} or $(\bf H_0)$  holds.

\end{rem}
Next, we discuss the number of eigenvalues  in the two following cases.

\medskip
\noindent
{\bf Case I: Suppose that ${\bf (H_0)}$ holds.} We can chose and fix two fundamental solutions $\phi$ and $\psi$ of $-y''+qy=0$ on $[0,1]$ satisfying that $\phi(0) = 0 = \psi(1)$, $\phi'(0) = 1$, $\psi'(1) = -1$, and  $W[\phi,\psi]=\omega \ne 0$. So, we have $\phi(1)=-\omega=\psi(0)$, and for any $\xi\in(0,1)$, 
\begin{equation}\label{eq:D-case1}
D_{1,\xi}=-\omega\psi(\xi),\quad D_{\xi,0}=-\omega\phi(\xi),
\end{equation}
which imply that ${\bf (H_1)}$ holds if and only if both of $\phi(t_i)$ and $\psi(t_i)$ are $0$ for any $i=1,2,\cdots,n$.

Let $\begin{pmatrix}
  \alpha_0\\
  \beta_0
\end{pmatrix} = \begin{pmatrix}
  1 \\
  0
\end{pmatrix}$. Because of \eqref{eq:alpha-K} and \eqref{eq:D-case1}, for any $1\leqslant K\leqslant n$, we have
\begin{equation}\label{eq:alpha-K2}
  \begin{aligned}
    \alpha_K =& \,\,
    1 + \frac{1}{\omega^2} \frac{\lambda}{\omega} \sum \limits_{k=1}^{K} m_k D_{t_{n+1},t_k} D_{t_k,t_0} + \frac{1}{\omega^2}\sum \limits_{l=2}^{K-1}  \left(\frac{\lambda}{\omega}\right)^l \sum \limits_{1 \le i_1 < \cdots < i_l \le K} \prod \limits_{k=1}^{l} m_{i_k} \prod \limits_{k=2}^{l} D_{t_{i_k}, t_{i_{k-1}}}  D_{t_{n+1}, t_{i_l}}  D_{t_{i_1}, t_0}  \\
    &+ \frac{1}{\omega^2}\left(\frac{\lambda}{\omega}\right)^K \prod \limits_{k=1}^K m_k \prod \limits_{k=1}^{K} D_{t_k, t_{k-1}} D_{t_{n+1}, t_K}.
  \end{aligned}
\end{equation}

In particular, when $K=n$, by \eqref{eq:alpha-K2}, we have that
\begin{equation}\label{eq:alpha}
  \begin{aligned}
    \alpha_n =& \,\,
    1 + \frac{1}{\omega^2} \frac{\lambda}{\omega} \sum \limits_{k=1}^{n} m_k D_{t_{n+1},t_k} D_{t_k,t_0} + \sum \limits_{l=2}^{n-1} \frac{1}{\omega^2} \left(\frac{\lambda}{\omega}\right)^l \sum \limits_{1 \le i_1 < \cdots < i_l \le n} \prod \limits_{k=1}^{l} m_{i_k} \prod \limits_{k=2}^{l} D_{t_{i_k}, t_{i_{k-1}}}  D_{t_{n+1}, t_{i_l}}  D_{t_{i_1}, t_0}  \\
    &+ \frac{1}{\omega^2}\left(\frac{\lambda}{\omega}\right)^n \prod \limits_{k=1}^n m_k \prod \limits_{k=1}^{n+1} D_{t_k, t_{k-1}}.
  \end{aligned}
\end{equation}
In addition, plugging $\psi(1)=0$ into \eqref{defn:charpoly}, we get the \emph{characteristics polynomial} $$p(\lambda) = \alpha_n\phi(1).$$

From \eqref{eq:alpha} and Theorem~\ref{thm:eigen-E-mat-poly}, we can deduce that
\begin{prop}\label{prop:case1}
$(a)$ If ${\bf(H_0)}$ holds but {\bf(H)} does not, then $0\leqslant \# \Spec {{\bf(E)}} < n$ and $0 \notin \Spec{\bf(E)}$.  Note that it is possible that $\# \Spec {{\bf(E)}} =0$, i.e., $\Spec {{\bf(E)}} =\emptyset$, for example,  Example~\ref{eg:emptyspec}(3).

\smnoind
$(b)$ If both of ${\bf(H_0)}$ and ${\bf(H)}$ hold, then $\# \Spec{{\bf(E)}}= n$, and $0 \notin \Spec{\bf(E)}$, i.e. {\bf(E)} has exactly $n$ non-zero distinct real eigenvalues.  
\end{prop}
\begin{prf}
The statement $(a)$ directly comes from the fact that if {\bf(H)} does not hold, then the factor $\prod \limits_{k=1}^{n+1} D_{t_k, t_{k-1}}$ in the coefficient of the term of highest degree in $\alpha_n$ is zero. 

For $(b)$, due to {\bf(H)}, we have that every $D_{t_k, t_{k-1}}$ $(k=1,2,\cdots,n+1)$ is non zero, and so the factor $\prod \limits_{k=1}^{n+1} D_{t_k, t_{k-1}}$ is non zero.

Set 
$$a_{k,k+1} = a_{k+1,k}:= - \frac{1}{D_{t_{k+1}, t_k}}, k = 1,2, \cdots, n, $$
and
$$a_{k,k} := \frac{D_{t_{k+1}, t_{k-1}}}{D_{t_{k+1},t_k}D_{t_k,t_{k-1}}}, k = 1, 2,\cdots, n.$$

By introducing the matrices 
\begin{equation}\label{eq:An-Mn}
A_n := (a_{i,j})_{n \times n} \quad\text{and}\quad  \tilde{M}_n := \begin{pmatrix}
  m_1 &        & O   \\
      & \ddots &     \\
  O   &        & m_n \\
\end{pmatrix},
\end{equation}
 where $a_{i,j}=0$ for other entries, we define
 \begin{equation}\label{eqn:alpha-K-det}
\tilde\alpha_K:=\prod \limits_{k=1}^{K+1} D_{t_k, t_{k-1}} \det\left(R_{[1,K]} \left(A_n + \frac{\lambda}{\omega}\tilde{M}_n\right)R_{[1,K]}\right), \quad \forall K\in\{1,2,\cdots,n\},
\end{equation}
where $R_{[1,K]}$ is the coordinate restriction to $\{1,2,\cdots,K\}$ of $\{1,2,\cdots,n\}$. Clearly, $R_{[1,n]}$ is the unit matrix.
 
By cofactor expansion, for any $1\leqslant K\leqslant n$, we have that
\begin{equation}\label{eqn:tildealpha-relation}
\tilde\alpha_K=\left(-\frac{D_{t_{K+1},t_K}}{D_{t_{K},t_{K-1}}}\right)\tilde\alpha_{K-2}+\left(\frac{D_{t_{K+1},t_{K-1}}}{D_{t_{K},t_{K-1}}}+\frac{\lambda}{\omega}m_KD_{t_{K+1},t_{K}}\right)\tilde\alpha_{K-1}.
\end{equation}
  
Then, applying \eqref{eq:alpha-K2} and Lemma~\ref{lem:Discrim-relation}, we can check, by induction, that, for $1\leqslant K\leqslant n$,  
\begin{equation}\label{eqn:tildealpha-K}
\begin{aligned}
    \tilde\alpha_K =& \,\,
    D_{t_{K+1},t_0} + \frac{\lambda}{\omega} \sum \limits_{k=1}^{K} m_k D_{t_{K+1},t_k} D_{t_k,t_0} +\sum \limits_{l=2}^{K-1}  \left(\frac{\lambda}{\omega}\right)^l \sum \limits_{1 \le i_1 < \cdots < i_l \le K} \prod \limits_{k=1}^{l} m_{i_k} \prod \limits_{k=2}^{l} D_{t_{i_k}, t_{i_{k-1}}}  D_{t_{K+1}, t_{i_l}}  D_{t_{i_1}, t_0}  \\
    &+ \left(\frac{\lambda}{\omega}\right)^K \prod \limits_{k=1}^K m_k \prod \limits_{k=1}^{K+1} D_{t_k, t_{k-1}}.
  \end{aligned}
\end{equation}

Consequently, since $D_{t_{n+1},t_0}=D_{1,0}=\omega^2$ here, when $K=n$, combined with \eqref{eq:alpha} and \eqref{eqn:tildealpha-K}, the equation \eqref{eqn:alpha-K-det} implies that
\begin{equation*}
\omega^2\alpha_n=\tilde\alpha_n=\prod \limits_{k=1}^{n+1} D_{t_k, t_{k-1}} \det\left(A_n + \frac{\lambda}{\omega}\tilde{M}_n\right),
\end{equation*}
and hence,  it follows that
\begin{equation}\label{eq:alpha-charpoly}
p(\lambda)=\phi(1)\alpha_n = -\frac{1}{\omega} \prod \limits_{k=1}^{n+1} D_{t_k, t_{k-1}} \det\left(A_n + \frac{\lambda}{\omega}\tilde{M}_n\right)=\frac{(-1)^{n+1}}{\omega^{n+1}} \prod \limits_{i=1}^{n}m_i \prod \limits_{k=1}^{n+1} D_{t_k, t_{k-1}} \det(X -\lambda I), 
\end{equation}
where $I$ is a unit, and
$$X=\tilde M_n^{-\frac{1}{2}}(-\omega A_n)\tilde M_n^{-\frac{1}{2}}.$$

Then, applying Theorem~\ref{thm:eigen-E-mat-poly} and Lemma~\ref{lem:dim-eigenspc}, the statement $(b)$ is proved, since the matrix $X$ is symmetric. 
\end{prf}

\begin{rem}\label{rem:alpha-charmat}
From the proof of Proposition~\ref{prop:case1}, we see that, if both ${\bf(H_0)}$ and ${\bf(H)}$ hold, we can see from \eqref{eq:alpha-charpoly} that the symmetric matrix $$X=\tilde M_n^{-\frac{1}{2}}(-\omega A_n)\tilde M_n^{-\frac{1}{2}}$$ is the corresponding \emph{characteristics matrix} of {\bf (E)}.
\end{rem}

\begin{eg}\label{eg:emptyspec}  
Let $\phi(x)=\frac{2}{3\pi}\sin \frac{3\pi}{2}x$ and $\psi(x)=-\frac{2}{3\pi}\cos \frac{3\pi}{2}x$ be the two linearly independent solutions of $-y''(x)-\frac{9\pi^2}{4}y(x) = 0$  on $[0,1]$. 
Their Wronskian determinant $W[\phi,\psi]$ is $\frac{2}{3\pi}$.
Then we consider the following three Dirichlet problems.

\smnoind
$(1)$ ${\bf (E_1):}$ $-y''(x)-\frac{9\pi^2}{4}y(x)=\lambda[\delta(x-\frac{1}{3}) + \delta(x - \frac{2}{3})]y(x)$ on $[0,1]$, $y(0) = y(1) = 0$. Here $\vec{t}=(t_0,t_1,t_2,t_3)=(0,\frac{1}{3},\frac{2}{3},1)$ and $\vec{m}=(m_1,m_2)=(1,1)$. Some simply calculation shows that $D_{t_3, t_0}=\frac{4}{9\pi^2}$, $D_{t_3, t_1}=D_{t_2, t_0}=0$, $D_{t_3, t_2}=D_{t_2, t_1}=D_{t_1, t_0}=-\frac{4}{9\pi^2}$, and so $p(\lambda)=-\frac{2}{3\pi}(1-\frac{4}{9\pi^2}\lambda^2)$. Consequently, we can see that, both of ${\bf (H)}$ and ${\bf(H_0)}$ hold, and $\Spec{{\bf (E_1)}} = \{\pm\frac{3\pi}{2}\}$, i.e., $\# \Spec{{\bf (E_1)}} = 2$.

\smnoind
$(2)$ ${\bf (E_2):}$ $-y''(x)-\frac{9\pi^2}{4}y(x) = \lambda[\delta(x-\frac{1}{4}) + \delta(x - \frac{1}{3})]y(x)$, $y(0) = y(1) = 0$. We can know that $(\bf H_0)$ holds but ${\bf (H)}$ is not true, because $D_{1,0}=\frac{4}{9\pi^2}$ and $D_{1,\frac{1}{3}}=0$. Furthermore, following the same procedures as in $(1)$, we can get $p(\lambda)=-\frac{2}{3\pi}(1-\frac{\lambda}{3\sqrt{2}\pi})$, and so $\Spec{{\bf (E_2)}} = \{3\sqrt{2} \pi \}$, i.e., $\# \Spec{{\bf (E_2)}} = 1<2$.

\smnoind
$(3)$ ${\bf (E_3):}$ $-y''(x)-\frac{9\pi^2}{4}y(x) = \lambda \delta(x - \frac{2}{3})y(x)$, $y(0) = y(1) = 0$. For this problem, we also can find that ${\bf (H_0)}$ hold but ${\bf (H)}$ fails, since $D_{1,0}=\frac{4}{9\pi^2}$ and $D_{\frac{2}{3},0}=0$. Hence, we have $p(\lambda)\equiv -\frac{2}{3\pi}$, which implies an extreme result that $\Spec{{\bf (E_3)}} = \emptyset$, i.e., $\# \Spec{{\bf (E_3)}} = 0$.
\end{eg}

\bigskip
\noindent
{\bf Case II: Suppose that ${\bf (H_0)}$ does not hold.}  
Since ${\bf (H_0)}$ does not hold, we can find non-zero number $\omega$ and a solution $\phi$ of $-y''+qy=0$ on $[0,1]$ such that $\phi(0) = 0=\phi(1)$, $\phi'(0)=1$ and $\phi'(1)=-\omega\ne 0$. 
Then we extend $\phi$ to a base $\{\phi,\ \psi\}$ of solution space of $-y''+qy=0$ such that $\psi(1)=1$. So, it follows that $W[\phi, \psi]=\omega$ and $\psi(0) =-\omega$. 
Then, we obtain the \emph{characteristics polynomial} $p(\lambda)=\beta_n $.

Moreover, for any $\xi\in(0,1)$, it is easily checked that
\begin{equation}\label{eq:D-case2}
D_{1,\xi}=-\phi(\xi),\quad D_{\xi,0}=-\omega\phi(\xi),
\end{equation}
which imply that, ${\bf(H_1)}$ holds, if and only if $\phi$ given above equals to $0$ at any nodal point $t_i$ $(\forall\ i=1,\cdots,n)$, if and only if the restriction of $\phi$ on $[t_i,t_j]$ is a non-zero solution of ${\bf (E_{[t_i,t_j]})}$, equivalently, $D_{t_j,t_i}=0$, for any $0\leqslant i<j\leqslant n+1$.

Set $\begin{pmatrix}
  \alpha_0 \\
  \beta_0 \\
\end{pmatrix} = \begin{pmatrix}
  1 \\
  0 \\
\end{pmatrix}$.  By \eqref{eq:beta-K} and \eqref{eq:D-case2}, we have that

\begin{equation}\label{eq:beta}
  \beta_n = -\frac{1}{\omega} \left[\begin{aligned}\frac{\lambda}{\omega} \sum_{k=1}^{n} m_k D_{t_{n+1}, t_k} D_{t_k, t_0} &+ \sum_{l=2}^{n-1}(\frac{\lambda}{\omega})^l \sum_{1 \le i_1 < \cdots < i_l \le n} \prod_{k=1}^{l}m_{i_k} \prod_{k=2}^{l} D_{t_{i_k}, t_{i_{k-1}}}  D_{t_{n+1}, t_{i_l}}  D_{t_{i_1}, t_0} \\
  &+ (\frac{\lambda}{\omega})^n \prod_{k=1}^{n}m_k \prod_{k=1}^{n+1}D_{t_{k}, t_{k-1}} \end{aligned}\right],
\end{equation}
which implies the following result.

\begin{prop}\label{prop:case2} 

\smnoind
$(a)$ Assume that neither $(\bf H)$ nor $(\bf H_0)$ holds. Then, when $(\bf H_1)$ holds, we have $\Spec{{\bf(E)}}=\BC$ (see Example~\ref{eg:C-spec}); when $(\bf H_1)$ does not hold, we have $1\leqslant \# \Spec{{\bf(E)}} < n$ and $0 \in \Spec{{\bf(E)}}$. 

\smnoind
$(b)$ If $(\bf H_0)$ does not hold but $(\bf H)$ holds, then $\# \Spec{{\bf(E)}}= n$ and $0 \in \Spec{{\bf(E)}}$, i.e. $\Spec{{\bf(E)}}= \{0\} \cup \{n-1 \text{ non-zero different real eigenvalues}\}$. 
\end{prop}
\begin{prf} Since $(\bf H_0)$ does not hold, the problem ${\bf (E_0)}$ has a non-zero solution, i.e., the problem ${\bf (E)}$ must have zero eigenvalue.

First, if $(\bf H_1)$ holds, then $p(\lambda)=\beta_n\equiv 0$, that is, every complex number is a root of $p(\lambda)=0$, and so $\Spec{{\bf(E)}}=\BC$. 

Next, we assume that $(\bf H_1)$ never holds. 
Then due to \eqref{eq:beta-K} and \eqref{eq:D-case2}, the term of degree $1$ in $\beta_n$ is non-zero unless $\lambda=0$, and so $p(\lambda)=\beta_n$ must be a non-zero polynomial.
The reason of the remainder of $(a)$ is the same as that in the proof of Proposition~\ref{prop:case1}(a).

For $(b)$, due to {\bf(H)}, we have that every $D_{t_k, t_{k-1}}$ $(k=1,2,\cdots,n+1)$ is non zero, and so the factor $\prod \limits_{k=1}^{n+1} D_{t_k, t_{k-1}}$ is also non zero.

Set $$b_{k,k}: = \frac{D_{t_{k+1},t_{k-1}}}{D_{t_{k+1},t_k}D_{t_k,t_{k-1}}}, k = 1, 2, \cdots, n, $$ and 
$$b_{k,k+1} = b_{k+1,k}:=-\frac{1}{D_{t_{k+1},t_k}}, k = 1, 2, \cdots, n.$$

Let $\tilde M_n$ be defined as in the proof of Proposition~\ref{prop:case1}(b), and $B_n$ be the matrix $(b_{i,j})_{n \times n}$ where $b_{i,j}=0$ for other entries. 
Here although $B_n$ has the same form as $A_n$ defined in the proof of Proposition~\ref{prop:case1}(b),  they are not the same one, because the fundamental solutions $\phi$ and $\psi$ here, which determine these discriminants $D_{t_{k+1},t_{k-1}}$ and $D_{t_{k+1},t_k}$ above, are both different from those in {\bf Case I}.

Then, noticing that $D_{t_{n+1},t_0}=D_{1,0}=0$, we can verify, by the same argument in the proof of Proposition~\ref{prop:case1}(b), that 
\begin{equation}\label{eq:beta-charpoly}
p(\lambda)=\beta_n = -\frac{1}{\omega}\prod \limits_{k=1}^{n+1}D_{t_k, t_{k-1}} \det\left(B_n + \frac{\lambda}{\omega} \tilde{M_n}\right)=\frac{(-1)^{n+1}}{\omega^{n+1}} \prod \limits_{i=1}^{n}m_i  \prod \limits_{k=1}^{n+1} D_{t_k, t_{k-1}} \det(X -\lambda I), 
\end{equation}
where $I$ is a unit matrix, and
$$X=\tilde M_n^{-\frac{1}{2}}(-\omega B_n)\tilde M_n^{-\frac{1}{2}}.$$
It is apparent that $X(\lambda)$ is the corresponding symmetric \emph{characteristics matrix} of {\bf (E)}.

Then, the rest of the proof of the statement $(b)$ follows from the same lines of argument as that of Proposition~\ref{prop:case1}(b).
\end{prf}

\begin{eg}\label{eg:0spec} 
Let $\phi(x)=\frac{1}{\pi}\sin\pi x$ and $\psi(x)=-\cos\pi x$ be the two linearly independent solutions of $-y''(x)-\pi^2y(x) = 0$  on $[0,1]$. 
Their Wronskian determinant $W[\phi,\psi]$ is $1$.
Then, following the similar arguments in Example~\ref{eg:emptyspec}(1), we can work out the eigenvalues of the following two Dirichlet problems, in which $(\bf H)$ holds while $(\bf H_0)$ does not hold.

\smnoind
$(1)$ For ${\bf (E_4):}$ $-y''(x)-\pi^2y(x) = \lambda \delta(x-\frac{1}{2})y(x)$ on $[0,1]$, $y(0) = y(1) = 0$, we have $p(\lambda)=-\frac{\lambda}{\pi^2}$, and $\Spec{{\bf (E_4)}} = \{0\}$, i.e., $\# \Spec{{\bf (E_4)}} = 1$.

\smnoind
$(2)$ For ${\bf (E_5):}$ $-y''(x)-\pi^2y(x) = \lambda(\delta(x-\frac{1}{4}) + \delta(x-\frac{1}{2}))y(x)$  on $[0,1]$, $y(0) = y(1) = 0$, we have $p(\lambda)=\frac{1}{2\pi^3}(\lambda^2-3\pi\lambda)$, and $\Spec{{\bf (E_5)}} = \{0, 3\pi\}$, i.e., $\#\Spec{{\bf (E_5)}} = 2$.
\end{eg}

\begin{lem}\label{lem:nonrealeigenvalue}
If $\lambda$ is a non-real eigenvalue of {\bf (E)} and $E(x)$ is an eigenfunction of $\lambda$, then $E(t_i)=0$ for any $i=0,1,2,\cdots,n,n+1$.
\end{lem}
\begin{prf}
Since the potential $q(x)$ is real value and $\lambda$ is a non-real eigenvalue of {\bf (E)}, it follows from the standard proof (cf. the proofs of \cite[Theorem 4.1.1]{Zet05} and \cite[Lemma 3.1]{GQ20}) that $\int_0^1 \delta_{n,\vec{m},\vec{t}}(x)\abs{E(x)}^2dx=0$. Hence $E(t_i)=0$ $(\forall\,i=0,1,2,\cdots,n,n+1)$, because of the boundary condition \eqref{defn:Eb} and $\int_0^1 \delta_{n,\vec{m},\vec{t}}(x)\abs{E(x)}^2dx=\sum^n_{i=1}m_i\abs{E(t_i)}^2$.
\end{prf}

In conclusion, we give the following complete result on numbers of Dirichlet eigenvalues and the existence of non-real eigenvalues.
\begin{thm}\label{thm:num-eigen} 
$(a)$ If the weight $w$ is $\delta_{n,\vec{m},\vec{t}}$ defined in {\bf(E)}, then either the number of its Dirichlet eigenvalues is at most $n$, or every complex number is a Dirichlet eigenvalue. More precisely, one and only one of the possibilities in both of Proposition~\ref{prop:case1} and Proposition~\ref{prop:case2} occurs.

\smnoind
$(b)$ Conversely, if the number of Dirichlet eigenvalues of the problem ${\bf(E_D)}$ is $K<+\infty$, then the weight $w(x)$ has the form of $\delta_{n,\vec{m},\vec{t}},$ and $n\geqslant K$.

\noindent
$(c)$ Suppose that {\bf(H)} holds.  Then the number of Dirichlet eigenvalues of the problem ${\bf(E_D)}$ is $K$, if and only if the weight $w(x)$ has the form of $\delta_{n,\vec{m},\vec{t}},$ and $n=K$. 

\noindent
$(d)$ Assume that both of {\bf(H)} and ${\bf(H_0)}$ hold. Then ${\bf(E_D)}$ has $K(<+\infty)$ distinct non-zero Dirichlet eigenvalues, if and only if the weight $w(x)$ has the form of $\delta_{n,\vec{m},\vec{t}},$ and $n=K$. 

\noindent
$(e)$ The problem {\bf (E)} has a non-real eigenvalue if and only if $\Spec({\bf E})=\BC$, if and only if  ${\bf (H_1)}$ holds. It means that $\Spec({\bf E})\subseteq\BR$ in the other cases.
\end{thm}
\begin{prf}
All of the statements $(a)$, $(c)$ and $(d)$ are the direct consequences of Proposition~\ref{prop:case1} and Proposition~\ref{prop:case2}. 

For $(b)$, by Remark~\ref{rem:mes-form}, the equation \eqref{defn:E-D} in ${\bf(E_D)}$ can be viewed as a measure differential equation. 
The statement $(b)$ is a generalization of the converse of \cite[Theorem 3.5(i)]{ZWMQX18} to the case that the potential is non-zero. In fact, the converse of \cite[Theorem 3.5(i)]{ZWMQX18} still holds even though the potential is non-zero.
Indeed, if we directly add a non-zero potential $q$ into the proof of the most crucial lemma \cite[Lemma 3.4]{ZWMQX18} for proving the converse of \cite[Theorem 3.5(i)]{ZWMQX18}, then
it can be easily seen, from the calculations in the proof of \cite[Lemma 3.4]{ZWMQX18}, that this added $q$ only generates an inconsequential infinitesimal $\frac{\norm{q}_{L^1}}{\sqrt{\tau}}$ (here $\tau$ is a sufficiently large number) in the key step  \cite[(3.25)]{ZWMQX18}, and the rest of the proof has no change. Hence, we can directly prove $(b)$ by imitating the proof of \cite[Theorem 3.5(ii)]{ZWMQX18}.

For $(e)$, by Remark~\ref{lem:nonrealeigenvalue}, we see that if $\lambda$ is a non-real eigenvalue of {\bf (E)}, then the corresponding eigenfunction $E(x)$ is zero at any $t_i$ $(\forall\, i=0,1,2,\cdots,n,n+1)$. So this $E(x)$ is a non zero solution of ${\bf (E_0)}$, and the restriction $E(x)|_{[t_i,t_j]}$ of $E(x)$ onto $[t_i,t_j]$ $(\forall\, 0\leqslant i<j\leqslant n+1)$ is a non zero solution of ${\bf (E_{[t_i, t_j]})}$. Hence, it is easily seen that ${\bf (H_1)}$ holds by Remark~\ref{rem:relation-hypo}, and then $\Spec({\bf E})=\BC$ by Proposition~\ref{prop:case2}(a).
\end{prf}

\begin{rem}\label{rem:general-result}
$(a)$ Compared with \cite[Problem (1.5)]{ZWMQX18}, the following more general problem is considered now:
\begin{equation}\label{eqn:gener-mes-form}
{\bf(E_M):}\quad\quad -dy^{\bullet}(x)+y(x)\,d\mu(x)=\lambda y(x)\,d\omega(x), \quad y(0)=y(1)=0,
\end{equation}
where $y(x)$ is defined as in \cite[(1.3)-(1.4)]{ZWMQX18}, $y^{\bullet}(x)$ stands for the generalized right derivative of $y$, $\lambda$ is the spectral parameter, $\mu(x)$ is an arbitrary absolutely continuous measure on $[0,1]$ with respect to Lesbegue measure $dx$, and $\omega(x)$ is an arbitrary real Radon measure on $[0,1]$.
Obviously, our problem ${\bf(E_D)}$ is an example of ${\bf(E_M)}$ without singularly continuous part.
 
From the explanation in the proof of Theorem~\ref{thm:num-eigen}(b), we know that \cite[Theorem 3.5(ii)]{ZWMQX18} also holds for ${\bf(E_M)}$, while \cite[Theorem 3.5(i)]{ZWMQX18}, which is no longer true for ${\bf(E_M)}$, has been updated with Theorem~\ref{thm:num-eigen} in the present paper.

\smnoind
$(b)$ From the proofs of Propositions~\ref{prop:case1}(b) and \ref{prop:case2}(b), we can see that, under the hypotheses {\bf(H)}, solving $p(\lambda)=0$ is equivalent to  solving $\det(X -\lambda I)=0$. 
Hence, due to the structure of $X$, we can simplify the 2nd step of the algorithm in Remark~\ref{rem:alg}, namely, we only need to calculate $2n$ discriminants: $D_{t_{k+1},t_{k-1}}$ and $D_{t_{k+1},t_k}$ ($k=1,2,\cdots,n$). 
\end{rem}

\medskip
\subsection{Two special cases.}
In this section, we discuss two simple but important cases.

\subsubsection{The case that the potential of {\bf(E)} is a zero function}\label{subsubsec:zero-q}
Consider
\begin{equation}\label{defn:Ed}
 {\bf(E_d):}\quad\quad  -y'' = \lambda \Delta_{n,\vec{m},\vec{t}}y,\quad y=y(x)\text{ on }[0,1],\quad y(0) = y(1) = 0,
 \end{equation}
which is a special case of {\bf(E)} with $q=0$. This equation~\eqref{defn:Ed} was also given in the setting of MDE in \cite[(3.13)]{ZWMQX18}. Solutions of \eqref{defn:Ed} are piecewise linear functions and its eigenvalue problem ${\bf(E_d)}$ is reduced to a linear system in $\BR^n$. So, in \cite[Lemma 3.3]{ZWMQX18}, the authors, using the difference method, equated ${\bf(E_d)}$ with a matrix eigenvalue problem. But, the difference method is not essential, and does not work in the general case. Next, we are going to use our approach to deal with ${\bf(E_d)}$.

For the homogenous equation $y''=0$ corresponding to ${\bf(E_d)}$, we can find two  fundamental solutions $\phi(x) = x$ and $\psi(x) = 1 - x$ satisfying $\omega:=W[\phi,\psi] = -1$, $\phi(0) = 0 = \psi(1)$, $\phi'(0) = 1$ and $\psi'(1) = -1$. More importantly, for any $0 \leqslant\xi < \eta \leqslant 1$, one has $D_{\eta, \xi} = \eta - \xi \ne 0$, which means that both of {\bf(H)} and ${\bf(H_0)}$ automatically hold for ${\bf(E_d)}$. 
Hence, it follows from Proposition~\ref{prop:case1}(b) that ${\bf(E_d)}$  must have exactly $n$ non-zero different eigenvalues.

Moreover, we can obtain that, $D_{t_k, t_{k-1}} = t_k - t_{k-1} \ne 0$ $(k=1,2,\cdots,n+1)$, and $D_{t_{k+1},t_{k-1}} = t_{k+1} - t_{k-1} \ne 0$ $(k=1,2,\cdots,n)$, and so, by \eqref{eq:alpha-charpoly}, we have
$$p(\lambda)=\alpha_n=\prod \limits_{i=1}^{n}m_i \prod_{k=1}^{n+1}(t_k - t_{k-1})\det(X - \lambda I),$$
where where $I$ is a unit, and
$X=\tilde M_n^{-\frac{1}{2}}A_n\tilde M_n^{-\frac{1}{2}}.$ Here $A_n$ and $\tilde M_n$ are both defined as in \eqref{eq:An-Mn}.  
 
Recalling \cite[(3.14)-(3.17)]{ZWMQX18}, the diagonal matrix $\tilde{M_n}$ above is denoted by $R$ in \cite{ZWMQX18}, the symmetric matrix $A_n$ above happens to be $A$ in \cite{ZWMQX18}, and $\hat A$ in \cite{ZWMQX18} is precisely the characteristics matrix $X$ here.

Therefore, we can see that our method and the difference method coincides in this special case that $q=0$. Differently, our method is useful for more general cases, and offer a general algorithm for finding eigenvalues.

\medskip
\subsubsection{The case that the weight is  of {\bf(E)} is a singel Dirac distribution}\label{subsubsec:singelDirac}
For any fixed $t\in(0,1)$, consider the following special case of  {\bf(E)}:
  \begin{equation}\label{defn:E-singleDirac}
    {\bf (E_{t,q})}:\quad\quad -y''(x)+q(x)y(x) = \lambda \delta(x-t)y(x)\text{ on }[0,1],\quad y(0) = y(1) = 0,
  \end{equation}
where $q\in L^1([0,1],\BR)$.

Note that, for ${\bf (E_{t,q})}$, the hypothesis $(\bf H_1)$ holds if and only if neither {\bf (H)} or $(\bf H_0)$  holds.
So Propositions~\ref{prop:case1} and ~\ref{prop:case2} show that one and only one  of the following possibilities occurs:
  \begin{equation*}
    \left\{ \begin{array}{ll}
      (1^\circ).\ \Spec{{\bf (E_{t,q})}} = \{0\}, & \text{when }(\bf H) \text{ holds but }  (\bf H_0) \text{ does not}; \\
      (2^\circ).\ \Spec{{\bf (E_{t,q})}} = \BC, & \text{when neither }(\bf H_0) \text{ or }  (\bf H) \text{ holds}; \\
      (3^\circ).\ \Spec{{\bf (E_{t,q})}} = \emptyset, & \text{when }(\bf H_0) \text{ holds but }  (\bf H) \text{ does not};\\
      (4^\circ).\ \Spec{{\bf (E_{t,q})}} = \{ \text{a unique non-zero real eigenvalue} \}, & \text{when both } (\bf H_0) \text{ and } (\bf H) \text{ hold}.
    \end{array} \right.
  \end{equation*}

For the case $(2^\circ)$, we give a simple example at once. 

\begin{eg}\label{eg:C-spec} 
Let $\phi(x)=\frac{1}{2\pi}\sin2\pi x$ and $\psi(x)=\frac{2}{\pi}\cos 2\pi x$ be the two linearly independent solutions of $-y''(x)-4\pi^2y(x) = 0$  on $[0,1]$. 
Their Wronskian determinant $W[\phi,\psi]$ is $-\frac{2}{\pi}$.
For ${\bf (E_6):}$ $-y''(x)-4\pi^2y(x) = \lambda \delta(x-\frac{1}{2})y(x)$ on $[0,1]$, $y(0) = y(1) = 0$, we have  $D_{1,0}=D_{1,\frac{1}{2}}=D_{\frac{1}{2},0}=0$, and so $p(\lambda)\equiv 0$, that is, $\Spec{{\bf (E_6)}} = \BC$. And $\phi$ is a common eigenfunction corresponding to all complex numbers.
\end{eg}

But our most concern is the fourth case. For the case $(4^\circ)$, we can choose the two fundamental solutions $\phi$ and $\psi$ of $-y''+qy=0$ satisfying
\begin{equation}\label{eqn:initialcond-phi-psi}
\phi(0)=0=\psi(1),\ \phi'(0)=1,\ W[\phi,\psi]=\omega\ne 0.
\end{equation}
which implies $\psi(0)=-\omega$ and $\phi(1)\psi'(1)=\omega$.
Hence, we have
\begin{equation}\label{eq:D-singelDirac}
D_{1,\xi}=\phi(1)\psi(\xi),\quad D_{\xi,0}=-\omega\phi(\xi), \text{ for any } \xi\in(0,1).
\end{equation}

Consequently, because {\bf(H)} holds, by \eqref{eq:D-singelDirac}, we know 

\noindent
{\bf Fact:} \emph{$\phi(t)$ and $\psi(t)$ are both non zero for the above given point $t$. This means that if {\bf(H)} holds for all $t\in(0,1)$, then both of $\phi(t)$ and $\psi(t)$ don not have sign-change on $[0,1]$. In particular, $\phi>0$ on $(0,1)$.}

Combining \eqref{eq:alpha-K}-\eqref{defn:charpoly} with \eqref{eqn:initialcond-phi-psi}, we get the \emph{characteristics polynomial}
\begin{equation}\label{eq:charapoly-singleDirac}
p(\lambda)=\phi(1)\alpha_n=\phi(1)\left(1+\frac{\lambda}{\omega}\phi(t)\psi(t)\right),
\end{equation}
and hence, for the above given $t\in(0,1)$ and $q\in L^1([0,1],\BR)$, the \emph{unique Dirichlet eigenvalue} $\lambda(t,q)$ of ${\bf (E_{t,q})}$ is given, that is, 
\begin{equation}\label{eqn:spec-singleDirac}
\lambda(t,q)=\frac{-\omega}{\phi(t)\psi(t)}\ne 0.
\end{equation}

\begin{lem}\label{lem:eigen-pos}
If for all $t\in(0,1)$, both of  ${\bf (H_0)}$ and {\bf (H)}  hold for ${\bf (E_{t,q})}$, then $\lambda(t,q)>0$ for any $t\in(0,1)$. 
\end{lem}
\begin{prf} From the initial condition \eqref{eqn:initialcond-phi-psi} of $\phi$ and {\bf Fact} above, we know that $\phi>0$ identically on $(0,1)$. 
If $W[\phi,\psi]=\omega> 0$, then $\psi(0)=-\omega<0$, which, due to the above {\bf Fact} again, implies that $\psi<0$ on $(0,1)$. Hence, by \eqref{eqn:spec-singleDirac}, we have $\lambda(t,q)>0$ for any $t\in(0,1)$. 
If $W[\phi,\psi]=\omega<0$, then $\psi(0)=-\omega>0$. So, from the same reason as above, we know $\psi>0$ on $(0,1)$, and hence also $\lambda(t,q)>0$ for any $t\in(0,1)$.
\end{prf}

\bigskip
\section{An application to Inverse Dirichlet problems involving single Dirac distribution weights}\label{sec:inverse-optim}\label{sec:inverse-prob}
\medskip

From Section~\ref{subsubsec:singelDirac}, excluding the extreme case of infinite eigenvalues, we can find that, if there exists $t_0\in (0,1)$ satisfying $\lambda(t_0,q)=0$, then ${\bf (H_0)}$ must not hold, and so $\lambda(t)\equiv 0$ for all $t\in(0,1)$. 
Hence, if $\Spec {\bf (E_{t,q})}$ is neither empty or $\BC$ for all $t\in(0,1)$, then, for any $t\in(0,1)$, the corresponding eigenvalue $\lambda(t,q)$ either identically vanishes, or is never zero. In this section, we mainly concern the latter case.

Set $$\CP[0,1]:=\{q\in L^1([0,1],\BR):\ \text{both of } {\bf (H_0)}\text{ and } {\bf (H)} \text{ hold for }  {\bf (E_{t,q})},\ \forall t\in(0,1)\}.$$

From the equations \eqref{eqn:initialcond-phi-psi}-\eqref{eqn:spec-singleDirac}, we can see that, an integrable function $q$ belongs to $\CP[0,1]$, if and only if $\Spec {\bf (E_{t,q})}$ has a unique non-zero eigenvalue for any $t\in(0,1)$, if and only if $-y''(x)+q(x)y(x)=0$ has two linear independent solutions such that either of them has only one zero point at $x=0$ or $x=1$.
Note that $\CP[0,1]$ is a big set. For example, the zero function belongs to $\CP[0,1]$; if $q(x)\in L^1([0,1],\BR)$ and $q(x)>0$ for any $x\in[0,1]$, then $q$ is a element of $\CP[0,1]$, due to the oscillation theory for initial value problem (cf., \cite[Section 2.6]{Zet05} or \cite[Section 3.3]{Cao16}).

Now, we consider the following problem:
  \begin{equation}\label{defn:E-singleDirac_P}
    {\bf (E^{\CP}_{t,q})}:\quad\quad -y''(x)+q(x)y(x) = \lambda \delta(x-t)y(x)\text{ on }[0,1],\quad y(0) = y(1) = 0,
  \end{equation}
where $t\in(0,1)$ and $q\in\CP[0,1]$.

Let $\phi$ and $\psi$ be the solutions defined as in \eqref{eqn:initialcond-phi-psi}. 
Then, combined with the initial conditions (\ref{eqn:initialcond-phi-psi}) as well as Lemma~\ref{lem:eigen-pos},  the formula \eqref{eqn:spec-singleDirac} tells us the following lemma immediately.

\begin{lem}\label{lem:property-eigen} For the problem ${\bf (E^{\CP}_{t,q})}$, we have

\smnoind
$(a)$ $\lambda(t,q)$, with respect to $t\in(0,1)$, is continuous, and has second-order derivative. In particular, there exists a constant $C>0$ such that $\abs{\lambda'(t)}\leqslant C\lambda^2(t)$ for any $t\in(0,1)$.

\smnoind
$(b)$ $\lambda(t,q)$ tends to the infinity as $t$ goes to $0$ or $1$.

\smnoind
$(c)$ $\lambda(t,q)>0$ for any $t\in(0,1)$.
\end{lem}

\begin{defn}\label{def:spec-data}
Let $\lambda(t,q)$ be the unique eigenvalue of SLP ${\bf (E^{\CP}_{t,q})}$ $(0<t<1)$. We call $\{\lambda(t,q)\}_{t\in(0,1)}$ a set of {\bf\emph{complete spectral data}} of SLP ${\bf (E^{\CP}_{t,q})}$. If $\CI$ is a proper subset of $(0,1)$, the set of eigenvalues $\{\lambda(t,q)\}_{t\in\CI}$ is called a set of {\bf\emph{incomplete spectral data}} of SLP ${\bf (E^{\CP}_{t,q})}$ with respect to $\CI$. 
\end{defn}

For simplicity,  if there is no ambiguity, we omit $q$ in \eqref{eqn:spec-singleDirac} and write $\lambda(t,q)$ as $\lambda(t)$.

\begin{prop}\label{prop:recover-q}
For the problem ${\bf (E^{\CP}_{t,q})}$, a unique potential function $q$ can be reconstructed from a set of given complete spectral data  $\{\lambda(t)\}_{t\in(0,1)}$, and 
\begin{equation}\label{eqn:recover-q}
q(x)=-\frac{1}{2}\cdot\frac{\lambda''(x)}{\lambda(x)}+\frac{3}{4}\cdot\left(\frac{\lambda'(x)}{\lambda(x)}\right)^2+\frac{1}{4}\lambda^2(x) {\rm\quad  on\ }(0,1).
\end{equation}
\end{prop}
\begin{prf}
Let $\phi$ and $\psi$ be the fundamental solutions as defined as in \eqref{eqn:initialcond-phi-psi}, and set $\omega:=W[\phi,\psi]\ne 0$, which, combined with \eqref{eqn:spec-singleDirac}, implies that, for any $t,\ x\in(0,1)$, one has
\begin{equation}\label{eqn:spec-singleDirac1}
\lambda(t)=\frac{-\omega}{\phi(t)\psi(t)}, 
\end{equation}
and
\begin{equation}\label{eq:psi-by-phi}
\psi(x)=\phi(x)\int_x^1\frac{-\omega}{\phi^2(s)}\, ds.
\end{equation}

Combining two equations above, we have 
\begin{equation}\label{eqn:phi-lambda}
\phi^2(x)\int^1_x\frac{1}{\phi^2(s)}\, ds=\frac{1}{\lambda(x)},\quad \forall x\in(0,1).
\end{equation}
And, from {\bf Fact} in Section~\ref{subsubsec:singelDirac}, we also can see that both of $\phi$ and $\psi$ don not have sign-change, and are never zero on $(0,1)$.

Set $u(x)=\int^1_x\frac{1}{\phi^2(s)}\ ds$. Then, $u'(x)=-\frac{1}{\phi^2(x)}$, and hence
$$-\frac{u(x)}{u'(x)}=\frac{1}{\lambda(x)},$$
whose general solution is 
\begin{equation}\label{eqn:u}
u(x)=Ce^{-\int^x_a\lambda(s)\,ds},\quad\forall x\in(0,1),
\end{equation}
where $a$ is an arbitrarily given point in $(0,1)$, and $C$ is an arbitrary positive constant which must be strictly positive.
So, we have
$$\phi^2(x)=\frac{C'}{\lambda(x)}e^{\int^x_a\lambda(s)\,ds},$$
that is,
\begin{equation}\label{eqn:recover-phi}
\phi(x)=\frac{C''}{\sqrt{\lambda(x)}}e^{\frac{\int^x_a\lambda(s)\,ds}{2}},
\end{equation}
where $C'=\frac{1}{C}$ and $C''=C'^{\frac{1}{2}}$.  

Hence, since the $\phi$ above only depends on $\lambda(t)$, we can recover the desired potential
\begin{equation*}
q(x)=\frac{\phi''(x)}{\phi(x)},
\end{equation*}
which, by substitution of \eqref{eqn:recover-phi} for $q(x)$, equals to \eqref{eqn:recover-q} on $(0,1)$.
\end{prf}

\begin{rem}\label{rem:recover-q-by-psi}
$(a)$ From the proof of Proposition~\ref{prop:recover-q}, we can find that the recovered $q$ is independent of the choices of the constant $C$ and the point $a$, since both of the constant $C$ and the item $e^{\frac{\int^x_a\lambda(s)ds}{2}}$ finally will be eliminated.

\smnoind
$(b)$ In the proof of Proposition~\ref{prop:recover-q}, by \eqref{eqn:spec-singleDirac1} and \eqref{eqn:recover-phi}, we can get 
\begin{equation}\label{eqn:recover-psi}
\psi(x)=\frac{-\omega}{C''\sqrt{\lambda(x)}}e^{-\frac{\int^x_a\lambda(s)\,ds}{2}},
\end{equation}
where $a\in(0,1)$, and $C''$ is a constant defined as in \eqref{eqn:recover-phi}.
Then, we similarly can get the relation $\phi(x)=\psi(x)\int_0^x\frac{1}{\psi^2(t)}\ dt$, and by setting $u(x)=\int^x_0\frac{1}{\psi^2(s)}\ ds$ or directly plugging \eqref{eqn:recover-psi} into the formula below, we have
$$\frac{\psi''(x)}{\psi(x)}=\frac{\phi''(x)}{\phi(x)}=q(x) {\rm\quad  on\ }(0,1),$$
\end{rem}

Remark~\ref{rem:recover-q-by-psi}(b) and the proof of the above proposition implicitly include the following properties of $\lambda(t)$.

\begin{cor}\label{cor:eigen-lim}
For any given $a\in(0,1)$, one has $$0<\lim_{t\rightarrow 0^+}\frac{\lambda(t)}{e^{\int^a_t \lambda(s)\,ds}}<+\infty \quad \text{ and }\quad 0<\lim_{t\rightarrow 1^-}\frac{\lambda(t)}{e^{\int^t_a \lambda(s)\,ds}}<+\infty.$$ Moreover, one has that $\int^1_a \lambda(s)\,ds=+\infty$ and $\int^a_0 \lambda(s)\,ds=+\infty$. 
\end{cor}
\begin{prf}
Because of \eqref{eqn:recover-phi} and  \eqref{eqn:recover-psi}, the first part of this corollary follows from the fact that $0<\phi(1),\psi(0)<+\infty$.  The last sentence directly follows from the proved part and Lemma~\ref{lem:property-eigen}(b).
\end{prf}

\begin{defn}\label{defn:spec-like}
Let $f(t)$ be a continuous and real value function on $(0,1)$. The function $f(t)$ is called {\bf\emph{spectrum-like function}}, if it satisfies all properties listed as follows: 

\noindent
$(1)$ $f(t)$ has second-order derivative on $(0,1)$.

\noindent
$(2)$ $f(t)>0$ for any $t\in(0,1)$.

\noindent
$(3)$ $f(t)$ tends to the infinity as $t$ goes to $0$ or $1$, and $\frac{f'(t)}{f^2(t)}=O(1)$ $(t\rightarrow 0^+\text{ or }1^-)$.

\noindent
$(4)$ $0\ne\lim_{t\rightarrow 0^+}\frac{f(t)}{e^{\int^a_t f(s)\,ds}}<+\infty$ and $0\ne\lim_{t\rightarrow 1^-}\frac{f(t)}{e^{\int^t_a f(s)\,ds}}<+\infty$ for any $a\in(0,1)$. 
 
Denote by $\mathcal{SL}[0,1]$ the set of all spectrum-like functions on $(0,1)$.
\end{defn}

Note that Definition~\ref{defn:spec-like}(3) and (4) imply that $\int^1_a f(s)ds=+\infty$ and $\int^a_0 f(s)ds=+\infty$ for any $a\in(0,1)$.
Obviously, the full spectrum data $\lambda(t)$ on $(0,1)$ in Proposition~\ref{prop:recover-q} is an element in $\mathcal{SL}[0,1]$. The simplest example of a function in $\mathcal{SL}[0,1]$ is $\frac{1}{t(1-t)}$, which gives the complete spectrum information of $-y''(x)=\lambda\delta(x-t)y(x)\text{ on }[0,1], \ y(0)=y(1)=0$.

Further, we can find that a spectrum-like function is equivalent to a set of complete spectral data of SLP ${\bf (E^{\CP}_{t,q})}$ $(0<t<1)$  for some $q\in\CP[0,1]$.

\begin{thm}\label{thm:q-spec-like}
For any function $f(t)$ defined on $(0,1)$, we have that, the set $\{f(t)\}_{t\in(0,1)}$ is the complete spectral data of SLP ${\bf (E^{\CP}_{t,q})}$ for some $q\in\CP[0,1]$,  if and only if $f(t)$ is an element in $\mathcal{SL}[0,1]$. 

More precisely, for any given $f\in\mathcal{SL}(0,1)$, we can construct an integrable function $Q$ as follows:
\begin{equation}\label{eqn:f-Q}
Q(x)=-\frac{1}{2}\cdot\frac{f''(x)}{f(x)}+\frac{3}{4}\cdot\left(\frac{f'(x)}{f(x)}\right)^2+\frac{1}{4}f^2(x) {\rm\quad  on\ }(0,1),
\end{equation}
such that $\{f(t)\}_{t\in(0,1)}$ is a set of complete spectrum data of ${\bf (E^{\CP}_{t,Q})}$ $(0<t<1)$.
\end{thm}

\begin{prf}
The ``only if" part is a direct corollary of Lemma~\ref{lem:property-eigen}, Proposition~\ref{prop:recover-q} as well as Corollary~\ref{cor:eigen-lim}. 

Conversely, for any $f\in\mathcal{SL}[0,1]$, the function $Q$ is the one obtained by \eqref{eqn:f-Q}. 
Let $\{\lambda(t)\}_{t\in(0,1)}$ be the set of complete spectrum data of ${\bf (E^{\CP}_{t,Q})}$. 
For proving the ``if" part, it suffices to show that, for any $x\in(0,1)$, one has $\lambda(x)=f(x)$. 

By Definition~\ref{defn:spec-like}(4), let $C:=\lim_{t\rightarrow 0^+}f(t)^{-\frac{1}{2}}e^{\frac{\int^a_t f(s)\,ds}{2}}\in(0,+\infty)$.
For any fixed $a\in(0,1)$, set $\phi_0(x)=2Cf(x)^{-\frac{1}{2}}e^{\frac{\int^x_a f(s)\,ds}{2}}$ for any $x\in(0,1)$. It can be directly verified that $Q=\frac{\phi_0''}{\phi_0}$ on $(0,1)$. 
And Definition~\ref{defn:spec-like}(3)(4) immediately tell us that $\phi_0(0^+)=0$ and $0<\phi_0(1^-)<+\infty$.  So $\phi_0$ can be extended to a solution of  $-y''+Qy=0$ on $[0,1]$ satisfying $\phi_0(0)=0$.

Moreover,  the derivation of $\phi_0$ is 
\begin{equation}\label{eqn:phi-deriv}
\phi_0'(x)=-\frac{Cf'(x)}{f(x)^{\frac{3}{2}}e^{\frac{1}{2}\int^a_x f(s)\,ds}}+\frac{Cf(x)^{\frac{1}{2}}}{e^{\frac{1}{2}\int^a_x f(s)\,ds}}=-C\cdot\frac{f'(x)}{f^2(x)}\cdot\frac{f(x)^{\frac{1}{2}}}{e^{\frac{1}{2}\int^a_x f(s)\,ds}}+\frac{Cf(x)^{\frac{1}{2}}}{e^{\frac{1}{2}\int^a_x f(s)\,ds}}.
\end{equation}
Then, applying Definition~\ref{defn:spec-like}(3)(4) again, we can know that $\phi_0'(0^+)<+\infty$. 
In addition, the initial value of $\phi'_0$ at $x=0$ must be non-zero, since $\phi_0$ is a non-zero solution of $-y''+Qy=0$ on $[0,1]$ and $\phi_0(0)=0$. 
Hence we set $\phi:=\frac{\phi_0}{\phi'_0(0)}$, which is a solution of $-y''+Qy=0$ satisfying $\phi(0)=0$ and $\phi'(0)=1$.

Next, for any $x\in(0,1)$,  let $\psi(x)=\phi(x)\int_x^1\frac{1}{\phi^2(s)}\ ds,$  which means that 
$\left(\frac{\psi(x)}{\phi(x)}\right)'=-\frac{1}{\phi^2(x)}$, and so $\phi(x)\psi'(x)-\phi'(x)\psi(x)=-1.$
Hence, we have $W[\phi,\psi]=-1$, and $\psi(1^-)=0$. Also we can check that 
$$\psi''(x)+Q\psi(x)=\left(\phi(x)\int_x^1\frac{1}{\phi^2(s)}\, ds\right)''+\frac{\phi''(x)}{\phi(x)}\phi(x)=0.$$
From all discussed above, we have obtained two linearly independent solutions $\phi$ and $\psi$ of $-y''+Qy=0$ on $[0,1]$ satisfying the initial condition \eqref{eqn:initialcond-phi-psi}.  
Since Definition~\ref{defn:spec-like}(4) yields that
$$\int_x^1f(t)e^{-\int^t_a f(s)\,ds}\ dt=\int_x^1d\left(-e^{-\int^t_a f(s)\,ds}\right)=e^{-\int^x_a f(s)\,ds},$$
we deduce from \eqref{eqn:spec-singleDirac} that
$$\lambda(x)=\frac{-W[\phi,\psi]}{\phi(x)\psi(x)}=\frac{1}{\phi^2(x)\int_x^1\frac{1}{\phi^2(s)}\, ds}
=\frac{f(x)}{e^{\int^x_a f(s)\,ds}\int_x^1f(t)e^{-\int^t_a f(s)ds}\, dt}=f(x),$$
for any $x\in(0,1)$.
Now the proof is done.
\end{prf}

Theorem~\ref{thm:q-spec-like} has a number of distinctive and interesting physical applications, which will be partially illustrated by the following living example. 

\begin{eg}\label{eg:vibration} 
Recall the model \eqref{eqn:vib-mass}-\eqref{eqn:vib-mass-v} of a vibrating string at the beginning of Section~\ref{sec:pre}. Clearly, the problem ${\bf (E_{t,q})}$ in \eqref{defn:E-singleDirac} is exactly the model \eqref{eqn:vib-mass}-\eqref{eqn:vib-mass-v} by setting $p(x)=1$ and $m=1$.  
So, because of \eqref{eqn:spec-singleDirac}, the vibration system \eqref{eqn:vib-mass}-\eqref{eqn:vib-mass-v} equipped with some external potential $q$ in $\CP[0,1]$ has a unique non-zero frequency. 

Then, we have a natural question: for the system \eqref{eqn:vib-mass}-\eqref{eqn:vib-mass-v}, how do we obtain a vibration at an expected frequency $\omega_0$? 
Theorem~\ref{thm:q-spec-like} give us an alternative method,
that is, we may try to find a function $f\in\CP[0,1]$ such that $f(\xi)=\omega^2_0$, and then construct a desired potential $q$ by  \eqref{eqn:f-Q}. 
But, is such potential unique? If not, then can we pick up the one that has the smallest $L^1-$norm (minimal potential energy)?
Motivated by these above, we will further introduce some related questions in the next section.
\end{eg}

\bigskip
\section{Open questions}\label{sec:question}
\medskip

Proposition~\ref{prop:recover-q} and Theorem~\ref{thm:q-spec-like} tell us that, for the problem ${\bf (E^{\CP}_{t,q})}$ in \eqref{defn:E-singleDirac_P}, a unique potential is determined if and only if a set of complete spectral data is known. 
Nevertheless, we can not always obtain the full spectral information, in other words, usually only a set of incomplete spectral data is known. 
In this case, we cannot make sure that a potential is recovered uniquely, but may settle for second best. 
Thereupon,  we naturally will consider inverse problems of characterizing the set consisting of all eligible potentials reconstructed by given incomplete spectral data, and the optimization problems of calculating the infimums of the $L^1$-norm of all such potentials and finding the optimal elements attaining the infimums.

For any $q\in\CP[0,1]$, let $\lambda(t,q)$ be the unique eigenvalue of SLP ${\bf (E^{\CP}_{t,q})}$ $(0<t<1)$ in \eqref{defn:E-singleDirac_P}.  For any subset $\CU$ of $(0,1)$ and any set $\Lambda:=\{\mu_t\}_{t\in\CU}$ of real numbers, define
$$\Omega(\Lambda,\CU):=\{q\in L^1[0,1]\cap\CP[0,1]\ | \ \mu_t=\lambda(t,q),\ \forall t\in\CU\},$$
$$E(\Lambda,\CU):=\inf\{\norm{q}_{L^1}\ | \ q\in\Omega(\Lambda,\CU)\},$$
and
$$M(\Lambda,\CU):=\{q\in \Omega(\Lambda,\CU)\ | \ \norm{q}_{L^1}=E(\Lambda,\CU)\}.$$

In fact, what Proposition~\ref{prop:recover-q} and Theorem~\ref{thm:q-spec-like} say is that $\Omega(\Lambda,(0,1))$ is a singleton set for any $f\in\mathcal{SL}[0,1]$ and the corresponding set  $\Lambda=\{f(t)\}_{t\in(0,1)}$.

To end this paper, we introduce some questions, motivated by our work, which remains to be studied in the sequel.

\noindent
{\bf Question 1.}  Suppose that $\Lambda=\{\lambda(t)\}_{t\in\CI}$ is a set of incomplete spectral data  of SLP ${\bf (E^{\CA}_{t,q})}$ with respect to $\CI$.
What is $\Omega(\Lambda,\CU)$? How do we calculate $E(\Lambda,\CU)$? Can we find out $M(\Lambda,\CU)$? Here $\CU$ may be a finite subset, a infinitely discrete subset or a continuous subinterval of $(0,1)$.  

More generally, we may consider the following question.

\noindent
{\bf Question 2.}   Denote by $\Spec_q{{\bf(E)}}$ the set of all Dirichlet eigenvalues of the problem {\bf(E)} with the integrable potential $q$. Theorem~\ref{thm:num-eigen}(a) says that $\#\Spec_q{{\bf(E)}}\leqslant n$ under the hypothesis ${\bf (H_0)}$. Then,  for any set $\Lambda_k$ of $k$ different real numbers ($k\leqslant n$), how do we characterize 
$$\Omega(\Lambda_k):=\{q\in L^1[0,1]\ | \ \Lambda_k\subset\Spec_q{{\bf(E)}}\}?$$

If $\Omega(\Lambda_k)$ is not a singleton set, then can we give
$$E(\Lambda_k):=\inf\{\norm{q}_{L^1}\ | \ q\in\Omega(\Lambda_k)\}?$$
And if $E(\Lambda_k)$ exists, then what is
$$M(\Lambda_k):=\{q\in \Omega(\Lambda_k)\ | \ \norm{q}_{L^1}=E(\Lambda_k)\}?$$

\bigskip

\section*{Acknowledgement}
This research was supported by the NSF of China (Grants 11701327, 11771253, 11271229 and 11971262) and the NSF of Shandong Province (Grant ZR2019MA038). 
The first author would like to thank Mr. Xiao Hu (Capital Normal University, Beijing), Professor Bing Xie (Shandong University, Weihai), Mr. Shuxiang Ma (Shandong University, Weihai) and Dr. Liwei Yu (Tsinghua University, Beijing) for their generous help.

\end{document}